\documentclass{article}
\usepackage{amsmath, amsfonts, amssymb}
\usepackage{graphicx}
\usepackage{latexsym}
\usepackage{hyperref}
\usepackage{amsfonts}
\usepackage[all]{xy}

\newtheorem{theorem}{Theorem}[section]
\newtheorem{lemma}[theorem]{Lemma}

\numberwithin{equation}{section}

\newcommand{\be}{\begin{eqnarray}}
\newcommand{\en}{\end{eqnarray}}
\newcommand{\no}{\nonumber}

\begin{document}


\author{
{\bf\large Adriano Cavalcante Bezerra$^{*}$\emph{ and}  Changyu Xia$^{**}$} \\
{\it\small *Instituto Federal Goiano-Campus Trindade}\\
{\it\small 75389-269, Trindade - GO, Brazil}\\
{\it\small $^{**}$ Universidade de Bras\'{\i}lia, Departamento de Matem\'atica}\\
{\it\small 70910-900, Bras\'{\i}lia - DF, Brazil}\\
{\it\small e-mails: adriano.bezerra@ifgoiano.edu.br, xia@mat.unb.br}}
\title{\textbf{Sharp Estimates for the First Eigenvalues of the Bi-drifting Laplacian}}

\maketitle{}

\begin{abstract}
In the present paper we study some kinds of the problems for the bi-drifting Laplacian operator and get
some sharp lower bounds for the first eigenvalue for these eigenvalue problems
on compact manifolds with boundary (also called a smooth metric measure space) and weighted Ricci curvature bounded inferiorly.
\end{abstract}
\emph{keywords}: Drifiting Laplacian; Bakry-Emery Ricci Curvature; Eigenvalues

\section{Introduction}
\label{intro}
For a given complete $n$-dimensional Riemannian manifold
$(M, \langle,\rangle)$ with a metric $\langle,\rangle$, the triple $(M, \langle,\rangle, e^{-\phi}d\upsilon)$ is called a smooth metric measure space, where $\phi:M\rightarrow\mathbb{R}$ is a smooth real-valued function on $M$ and $d\upsilon$ is the Riemannian
volume element related to $\langle,\rangle$ (sometimes, we also call $d\upsilon$ the volume density). On
a smooth metric measure space $(M, \langle,\rangle, e^{-\phi}d\upsilon)$, we can define the so-called drifting
Laplacian (also called weighted Laplacian) $L_{\phi}$ as follows
\be\no
L_{\phi}:=\triangle-\langle \nabla \phi,\nabla (\cdot)\rangle,
\en
where $\nabla$ and $\triangle$ are the gradient operator and the Laplace operator, respectively.
Some interesting results concerning eigenvalues of the drifting Laplacian can be
found, for instance, in \cite{CMZ2,DB,DWLX,MD,XX}, among others.

On smooth metric measure spaces, we can also define the so-called $Bakry$-$Emery$ $Ricci$ tensor $Ric_\phi$ by
\be\label{1.1}
Ric_{\phi} = Ric + \nabla^{2}\phi,
\en
which is also called the weighted Ricci curvature. Here, $Ric$ is the Ricci curvature on $M$. The equation $Ric_\phi=k\langle,\rangle$ for some
constant $k$ is just the gradient Ricci soliton equation, which plays an important role
in the study of Ricci flow. For $k=0$, $k>0$ and $k<0$, the gradient Ricci soliton
$(M, \langle,\rangle, e^{-\phi}d\upsilon, k)$ is called steady, shrinking, and expanding, respectively.

In \cite{CCWX}, Chen-Cheng-Wang-Xia gave some lower bounds for the first eigenvalue of
four kinds of eigenvalue problems of the biharmonic operator on compact manifolds
with boundary and positive Ricci curvature, where two of them are in the direction of the buckling and clamped plate problems.
Posteriorly in \cite{DB}, Du and Bezerra extended the results of Chen-Cheng-Wang-Xia for the bi-drifting Laplacian operator. In particular, in the Theorems 1.7-1.9 they obtained lower estimates for the first eigenvalue of some eigenvalue problems for the bi-drifting Laplacian operator, defined in an smooth metric measure space $(M, \langle,\rangle, e^{-\phi} d\upsilon)$ with boundary $\partial M$ and a limiting condition in the Bakry-Emery Ricci curvature, if the weighted mean curvature $H_{\phi}$ of $\partial M$ is nonnegative. The weighted mean curvature $H_{\phi}$ will be defined in the next section.

In the first part of this paper, in the Theorem \ref{Theorem 1} and  Theorem \ref{Theorem 2} we will improve the results of Du and Bezerra, removing the condition at the weighted mean curvature of the boundary. The results are shown below.
\begin{theorem}\label{Theorem 1}
Let $(M, \langle,\rangle, e^{-\phi} d\upsilon)$ be an $n(\geq2)$-dimensional compact connected
smooth metric measure space with boundary $\partial M$ and denote by $\nu$ the outward unit
normal vector field of $\partial M$. Denote by $\lambda_{1}$ the first eigenvalue with Dirichlet boundary condition of the
drifting Laplacian of M and let $\Gamma_{1}$ be the first eigenvalue of the problem:
\be\label{1.2}
\left\{ \begin{array}{ll}
L^2_{\phi}u=\Gamma u \ \ \ \ \ in \ M,   & \\
u=\frac{\partial^{2} u}{\partial\nu^{2}}=0 \ \ on \ \partial M. & \\
\end{array} \right.
\en
Assume that
\be\label{1.3}
Ric_{\phi}\geq\dfrac{|\nabla \phi|^{2}}{na}+b,
\en
for some positive constants $a$ and $b$. Then we have
\be\label{1.4}
\Gamma_{1}>\lambda_{1}\left(\dfrac{\lambda_{1}}{n(a+1)}+b\right).
\en
\end{theorem}

Changing the first equation in \ref{1.2} and maintaining the boundary condition, we have  the following result:
\begin{theorem}\label{Theorem 2}
Under the assumption of Theorem \ref{Theorem 1}, let  $\eta_{1}$ the first
eigenvalue of the problem :
\be\label{1.5}
\left\{ \begin{array}{ll}
L^2_{\phi}u=-\eta L_{\phi} u \ \ in \ M,   & \\
u=\frac{\partial^{2} u}{\partial\nu^{2}}=0 \ \ \ \ on \ \partial M. & \\
\end{array} \right.
\en
Then we have
\be\label{1.6}
\eta_{1}>\dfrac{\lambda_{1}}{n(a+1)}+b.
\en
\end{theorem}

The eigenvalue problems (\ref{1.2}) and  (\ref{1.5})  should be compared with the clamped
plate problem and the buckling problem for drifting operator, respectively. The later two ones are as
follows:
\be\label{1.7}
\left\{ \begin{array}{ll}
L^2_{\phi}u=\Lambda u \ \ \ \ \ \ in \ M,   & \\
u=\frac{\partial u}{\partial\nu}=0 \ \ \ \ on \ \partial M. & \\
\end{array} \right.
\en
\be\label{1.8}
\left\{ \begin{array}{ll}
L^2_{\phi}u=-\eta L_{\phi} u \ \ in \ M,   & \\
u=\frac{\partial u}{\partial\nu}=0 \ \ \ \ \ on \ \partial M. & \\
\end{array} \right.
\en
The Theorem \ref{Theorem 1} and Theorem \ref{Theorem 2}  shows that the first eigenvalue of problems \ref{1.2}  and  \ref{1.5} are
closely related to the first Dirichlet eigenvalue of the drifting Laplacian. Next we will be interested in some types Steklov eigenvalue problems of the
bi-drifting Laplace operator. The eigenvalue problems we are interested in are as follows:
\be\label{1.9}
\left\{ \begin{array}{ll}
L^2_{\phi}u=0 \ in \ M,   & \\
u=L_{\phi}u-p\frac{\partial u}{\partial\nu}=0 \ on \ \partial M. & \\
\end{array} \right.
\en
\be\label{1.10}
\left\{ \begin{array}{ll}
L^2_{\phi}u=0 \ in \ M,   & \\
u=\frac{\partial^{2} u}{\partial\nu^{2}}-q\frac{\partial u}{\partial\nu}=0 \ on \ \partial M. & \\
\end{array} \right.
\en
\be\label{1.11}
\left\{ \begin{array}{ll}
L^2_{\phi}u=0 \ in \ M,   & \\
\frac{\partial u}{\partial \nu}=\frac{\partial (L_{\phi}u)}{\partial\nu}+\xi u=0 \ on \ \partial M. & \\
\end{array} \right.
\en
and
\be\label{1.12}
\left\{ \begin{array}{ll}
L^2_{\phi}u=0 \ in \ M,   & \\
\frac{\partial u}{\partial \nu}=\frac{\partial (L_{\phi}u)}{\partial\nu}+\beta \overline{\triangle}u+ \varsigma u=0 \ on \ \partial M. & \\
\end{array} \right.
\en

Elliptic problems with parameters in the boundary conditions are called Steklov
problems from their first appearance in \cite{ST}. The problem \ref{1.9} was considered by
Kuttler \cite{K} and Payne \cite{P} who studied the isoperimetric properties of the first
eigenvalue $p_{1}$ which is the sharp constant for $L^{2}$ a priori estimates for solutions
of the (second order) Laplace equation under nonhomogeneous Dirichlet boundary
conditions.  One can see that  $p_{1}$  is positive and given
by
\be\label{1.13}
p_{1} & = & \min_{w|_{\partial M}=0, \ w\neq const.} \frac{\int_{M} (L_{\phi}w)^{2}e^{-\phi}d\upsilon}
{\int_{\partial M} \left(\frac{\partial w}{\partial \nu}\right)^{2}e^{-\phi}dA}.
\en
The problem (\ref{1.10}) is a natural Steklov problem for the drifting Laplacian, and is equivalent to (\ref{1.9})  when
the mean curvature of $\partial M$ and $\phi$ are constants. In Theorem \ref{Theorem 3} we have a sharp relation between the first
eigenvalues of \ref{1.9} and \ref{1.10}.
\begin{theorem}\label{Theorem 3}
Let $(M, \langle,\rangle, e^{-\phi} d\upsilon)$ be an $n$-dimensional compact
smooth metric measure space with boundary $\partial M$ and non-negative Ricci Bakry-Emery curvature. Denote by $p_{1}$ and $q_{1}$ the first eigenvalue
of the problems (\ref{1.9}) and (\ref{1.10}), respectively. Then we have
\be\label{1.14}
q_{1}\geq\dfrac{p_{1}}{n(a+1)},
\en
with equality holding if and only if $M$ is isometric to a ball in $\mathbb{R}^{n}$.
\end{theorem}

The next result is a sharp lower bound for $p_{1}$.
\begin{theorem}\label{Theorem 4}
Let $(M, \langle,\rangle, e^{-\phi} d\upsilon)$ be an $n(\geq2)$-dimensional compact
smooth metric measure space. Suppose that
\be\label{1.15}
Ric_{\phi}\geq\dfrac{|\nabla \phi|^{2}}{na}-b,
\en
for some positive constants $a$ and $b$. Denote by $\lambda_{1}$ the first eigenvalue with Dirichlet boundary condition of the
drifting Laplacian of $M$ and let $p_{1}$ be the first eigenvalue of the
problem (\ref{1.9}). If the weighted mean curvature of $\partial M$ is bounded below by a
positive constant $c$, then we have
\be\label{1.16}
p_{1}\geq\dfrac{n(a+1)(n-1)c\lambda_{1}}{n(a+1)(\lambda_{1}+b)-\lambda_{1}},
\en
with equality holding if and only if $M$ is isometric to a Euclidean $n$-ball of radius $1/c$.
\end{theorem}

The problem (\ref{1.11}) was first studied in (\cite{KS}) where some estimates for the first
non-zero eigenvalue $\xi$ were obtained. When $M$ is an Euclidean ball, all the eigenvalues
of the problem (\ref{1.10}) have been recently obtained in \cite{QX2}. Also, the authors
proved an isoperimetric upper bound for $\xi_{1}$ when $M$ is a bounded domain in $\mathbb{R}^{n}$.
The Rayleigh-Ritz formula for $\xi_{1}$  is:
\be\label{1.17}
\xi_{1}= \ min \ \frac{\int_{M} (L_{\phi}w)^{2}e^{-\phi}d\upsilon}{\int_{\partial M} w^{2}e^{-\phi}dA},
\en
where $0\neq w\in H^{2}(M), \ \int_{\partial M}w=0=\partial_{\nu}w|_{\partial M}.$

The problem (\ref{1.12}) is a so called Wentzell problem for the bi-drifting laplace operator which
is motivated by (\ref{1.11}) and the following Wentzell-Laplace problem:
\be\label{1.18}
\left\{ \begin{array}{ll}
\triangle u=0 \  in \ M,   & \\
-\beta \overline{\triangle} u+\partial_{\nu} u=\lambda u \ \ on \ \partial M, & \\
\end{array} \right.
\en
where $\beta$ is a given non-negative number. The problem (\ref{1.18}) has been studied
recently, in \cite{DKL}, \cite{QX2}, etc.
The first non-zero eigenvalue of ((\ref{1.12})) can be characterized as
\be\label{1.19}
\zeta_{1,\beta}= \ min \ \frac{\int_{M} (L_{\phi}w)^{2}e^{-\phi}d\upsilon+\beta\int_{\partial M} |\overline{\nabla} w|^{2}e^{-\phi}dA}{\int_{\partial M} w^{2}e^{-\phi}dA},
\en
where $0\neq w\in H^{2}(M), \ \int_{\partial M}w=0=\partial_{\nu}w|_{\partial M}$.

From (\ref{1.17}), one can see that if $\beta>0$, $\xi_{1}$ is the first non-zero eigenvalue of the
Steklov problem (\ref{1.11}) and $\lambda_{1}$ the first non-zero eigenvalue of the drifting Laplacian of $\partial M$,
then we have
\be\label{1.20}
\zeta_{1,\beta}\geq\xi_{1}+\beta\lambda_{1},
\en
with equality holding if and only if any eigenfunction $f$ corresponding to $\zeta_{1,\beta}$ is an
eigenfunction corresponding to $\xi_{1}$ and $f|_{\partial M}$ is an eigenfunction corresponding to
$\lambda_{1}$. Our last result is a lower bound for $\zeta_{1}$.
\begin{theorem}\label{Theorem 5}  Let $(M, \langle,\rangle, e^{-\phi} d\upsilon)$ be an $n$-dimensional compact
smooth metric measure space with boundary $\partial M$ and suppose that (\ref{1.15}) is satisfied.  Assume that the
principal curvatures of $\partial M$ are bounded below by a positive constant $c$ and denote
by  $\zeta_{1}$ the first eigenvalue of the problem(\ref{1.12}). Then we have
\be\label{1.21}
\zeta_{1}>\dfrac{cn(a+1)\lambda_{1}\mu_{1}}{n(a+1)(\mu_{1}+b)-\mu_{1}}+\beta\lambda_{1},
\en
where $\mu_{1}$ and $\lambda_{1}$ are the first nonzero Neumann eigenvalue of the drifting Laplacian of $M$
and the first nonzero eigenvalue of the drifting Laplacian of $\partial M$, respectively.
\end{theorem}
\section{Proof of Theorems}
In this section, we will prove the theorems of the section 1. Before doing this, let us recall the Reilly$'$s formula. Let $M$ be an n-dimensional
compact manifold $M$ with boundary $\partial M$. We will often write $\langle,\rangle$ the Riemannian metric on $M$ as well as that
induced on $\partial M$. Let $\nabla$ and $\triangle$ be the connection and the Laplacian on $M$, respectively. Let $\nu$ be the unit outward normal
vector of $\partial M$. The shape operator of $\partial M$ is given by $S(X) =\nabla_{X}\nu$ and the second fundamental form of $\partial M$ is defined
as $II(X, Y ) = \langle S(X), Y \rangle$, here $X$, $Y\in\partial M$. The eigenvalues of S are called the principal curvatures of $\partial M$.

We will denote the weighted measure by $d\mu=e^{-\phi}dv$ and $d\vartheta=e^{-\phi}dA$ on $M$ and $\partial M$, respectively.  The weighted mean curvature as a natural generalization of the mean curvature for Riemann manifolds with density and is defined by
\be\label{2.1}
H_{\phi}=H-\frac{1}{n-1}\phi_{\nu},
\en
where $H$ denotes the usual mean curvature of $\partial M$, given by
$H = \dfrac{1}{n-1} tr S$, and $tr S$ denotes the trace of $S$.
For a smooth function $f$ defined on an n-dimensional
compact manifold $M$ with boundary, Ma and Du (\cite{MD}) extended the Reilly$'$s
formula for Riemann manifolds with density and showed that the following identity holds if
\be\label{2.2}
&  & \no\int_{M}(L_{\phi}f)^{2}-|\nabla^{2}f|^{2}-Ric_{\phi}(\nabla f,\nabla f)d\mu\\
& = &\int_{\partial M}2(\overline{L}_{\phi}f)f_{\nu}+(n-1)H_{\phi}(f_{\nu})^{2}+II(\overline{\nabla} f,\overline{\nabla }f)d\vartheta.
\en
Here $\overline{L}_{\phi}=\overline{\triangle}-\langle\overline{\nabla}\phi,\overline{\nabla}(.)\rangle$ and $|\overline{\nabla}^{2}f|$ are drifting operator and the Hessian of $f$ on $\partial M$, with respect to the induced metric
on $\partial M$, respectively.\\
\\
\textbf{\emph{Proof of Theorem} \ref{Theorem 1}.} Let $f$ be an eigenfuction of the problem (\ref{1.2}) corresponding to the first eigenvalue $\Gamma_{1}$, that is,
\be\label{2.3}
L^2_{\phi}f=\Gamma_{1}f \ \ in \ \ M, \ \ \ \
f=\frac{\partial^{2} f}{\partial\nu^{2}}=0 \ \ on \ \ \partial M.
\en
Multiplying the first equality in (\ref{2.3}) by $f$ and integrating on $M$, follows from the divergence theorem
\be\label{2.4}
\Gamma_{1}\int_{M}f^{2}d\mu=\int_{M}(L_{\phi}f)^{2}d\mu-\int_{\partial M} hL_{\phi}f d\vartheta,
\en
where $h=\dfrac{\partial f}{\partial \nu}|_{\partial M}$.
We now consider the following Lemma (\cite{BZ}):
\begin{lemma}\label{Lemma 1} Let $M$ an $n$-dimensional Riemannian manifold with boundary $\partial M$ and let $f\in C^{\infty}(M)$. Then for all $p\in M$ we have
\be\label{2.5}
L_{\phi}f = \overline{L}_{\phi}f+(n-1)H_{\phi}f_{\nu}+\nabla^{2}f (\nu,\nu),
\en
where $H_{\phi}$ is the weighted mean curvature of $\partial M$ and $L_{\phi}$ and $\overline{L}_{\phi}$ are the drifting Laplacian operators defined in $M$ and $\partial M$ respectively.
\end{lemma}
Since $f|_{\partial M}=\dfrac{\partial^{2} f }{\partial \nu^{2}}|_{\partial M}=0$ we have
\be\label{2.6}
L_{\phi}f\mid_{\partial M}=(n-1)H_{\phi}h.
\en
From Reilly's formula, we infer
\be\label{2.7}
\int_{M}(L_{\phi}f)^{2}-|\nabla^{2}f|^{2}-Ric_{\phi}(\nabla f,\nabla f) \ d\mu=(n-1)\int_{\partial M}H_{\phi}h^{2}d\vartheta.
\en
Combining (\ref{2.4}), (\ref{2.6}) and (\ref{1.3}), we get
\be\label{2.8}
\no\Gamma_{1}&=&\dfrac{\int_{M}(|\nabla^{2}f|^{2}+Ric_{\phi}(\nabla f, \nabla f)) \ d \mu}{\int_{M}f^{2}}\\
&\geq&\dfrac{\int_{M}(|\nabla^{2}f|^{2}+\left(\frac{\nabla \phi}{na}+b \right)|\nabla f|^{2}) \ d \mu}{\int_{M}f^{2}}.
\en
We get easily that
\be\label{2.9}
(\triangle f)^{2}=(L_{\phi}f+\langle\nabla\phi,\nabla f\rangle)^{2}\geq\dfrac{(L_{\phi}f)^{2}}{a+1}-\dfrac{\langle\nabla\phi,\nabla  f\rangle^{2}}{a}.
\en
The Schwarz inequality implies that
\be\label{2.10}
|\nabla^{2}f|^{2}\geq\dfrac{1}{n}(\triangle f)^{2},
\en
with equality holding if and only if
\be\label{2.11}
\nabla^{2}f=\dfrac{\triangle f}{n}\langle,\rangle. \
\en
It then follows from the Schwarz inequality that
\be\label{2.12}
|\nabla^{2}f|^{2}\geq\dfrac{1}{n}(\triangle f)^{2}\geq\dfrac{(L_{\phi}f)^{2}}{n(a+1)}-\dfrac{\langle\nabla\phi,\nabla  f\rangle^{2}}{na}
\en
Therefore, substituting (\ref{2.12}) in  (\ref{2.8}),  we have
\be\label{2.13}
\no\Gamma_{1}&\geq&\dfrac{\int_{M}\left(\dfrac{(L_{\phi}f)^{2}}{n(a+1)}-\dfrac{|\nabla\phi|^{2}|\nabla  f|^{2}}{na}+\left(\dfrac{|\nabla \phi|^{2}}{na}+b\right)|\nabla f|^{2}\right) d\mu}{\int_{M}f^{2}\ d\mu}\\
&=&\dfrac{\int_{M}\left(\dfrac{(L_{\phi}f)^{2}}{n(a+1)}+b|\nabla f|^{2}\right)\ d\mu}{\int_{M}f^{2}d\mu}.
\en
with equality holding if and only (\ref{2.11}) holds and
\be\no
Ric_{\phi}=\dfrac{|\nabla \phi|^{2}}{na}+b.
\en
On the other hand, since $f$ is not a zero function which vanishes on $\partial M$, we
know that
\be\label{2.14}
\int_{M}(L_{\phi} f)^{2}d\mu\geq\lambda_{1} \int_{M}|\nabla f|^{2}d\mu\geq\lambda_{1}^{2}\int_{M}f^{2}d\mu.
\en
with equality holding if and only if $f$ is a first eigenfunction of the Dirichlet problem for drifting Laplacian of $M$. Thus by (\ref{2.13}) and (\ref{2.14}) we conclude that
\be\label{2.15}
\Gamma_{1}\geq\lambda_{1}\left(\dfrac{\lambda_{1}}{n(a+1)}+b\right).
\en
suppose that $\Gamma_{1}=\lambda_{1}\left(\dfrac{\lambda_{1}}{n(a+1)}+b\right)$ is valid. Then (\ref{2.9}) becomes
\be\label{2.16}
(\triangle f)^{2}=(L_{\phi}f+\langle\nabla\phi,\nabla f\rangle)^{2}=\dfrac{(L_{\phi}f)^{2}}{a+1}-\dfrac{\langle\nabla\phi,\nabla  f\rangle^{2}}{a}.
\en
which means that  $\phi$ is not a constant and $\triangle f-\dfrac{1}{a}\langle\nabla f, \nabla \phi\rangle=0$
 holds everywhere on $M$. Multiplying the
above inequality by $f$ and integrating on $M$ with respect to $e^{\frac{1}{a}\phi} d\upsilon$ give that
\be\label{2.17}
0=\int_{M}f\left( \triangle f-\dfrac{1}{a}\langle\nabla f, \nabla \phi\rangle \right)e^{\frac{1}{a}\phi}d\upsilon=-\int_{M}|\nabla f|^{2}e^{\frac{1}{a}\phi}d\upsilon
\en
From above equality, we know that $f$ is a constant function on $M$, which is a contradiction since
$f$ is the first eigenfunction of bi-drifting Laplacian and cannot be a constant. Therefore, we have $\Gamma_{1}>\lambda_{1}\left(\dfrac{\lambda_{1}}{n(a+1)}+b\right)$.\hfill $\Box$ \vspace{0.2cm}
\\
\textbf{\emph{Proof of Theorem} \ref{Theorem 2}.} The discussions are similar to those in the proof of Theorem
\ref{Theorem 1}.  Let $g$ be the eigenfunction
of the problem (\ref{1.5}) corresponding to the first eigenvalue  $\eta_{1}$, that is,
\be\label{2.18}
L^2_{\phi}g=-\eta_{1}\triangle g \ \ in \ \ M, \ \ \ \
g=\frac{\partial^{2} g}{\partial\nu^{2}}=0 \ \ on \ \ \partial M.
\en
Multiplying the first equality in (\ref{2.18}) by $g$ and integrating on $M$, follows from the divergence theorem that
\be\label{2.19}
\eta_{1}\int_{M}|\nabla g|^{2}d\mu=\int_{M}(L_{\phi}g)^{2}d\mu-\int_{\partial M} sL_{\phi}g d\vartheta,
\en
where $s=\dfrac{\partial g}{\partial \nu}|_{\partial M}$. Also, we have
\be\label{2.20}
L_{\phi}g\mid_{\partial M}=(n-1)H_{\phi}s.
\en
Hence
\be
\no\eta_{1}=\dfrac{\int_{M}(L_{\phi}g)^{2}d\mu-(n-1)\int_{\partial M}H_{\phi} s^{2} d\mu}{\int_{M}|\nabla g|^{2}d\mu},
\en
which, by hyphotesis and by Reilly's formula and (\ref{2.12}) gives
\be\label{2.21}
\no\eta_{1}&\geq&\dfrac{\int_{M}\left(|\nabla^{2}g|^{2}+\left(\dfrac{|\nabla \phi|^{2}}{na}+b \right)|\nabla g|^{2}\right) \ d \mu}{\int_{M}|\nabla g|^{2}d\mu}\\
&\geq&\dfrac{\int_{M}\left(\dfrac{(L_{\phi}g)^{2}}{n(a+1)}+b|\nabla g|^{2}\right)\ d\mu}{\int_{M}|\nabla g|^{2}d\mu}.
\en
We can see that (\ref{2.14} also holds for $g$, and therefore
\be\no
\eta_{1}\geq\dfrac{\lambda_{1}}{n(a+1)}+b.
\en
In a similar way to what was done in the proof of the Theorem \ref{Theorem 1}, if we suppose that $\eta_{1}=\dfrac{\lambda_{1}}{n(a+1)}+b$, similarly to what was done,
\be\label{2.22}
0=\int_{M}g\left( \triangle g-\dfrac{1}{a}\langle\nabla g, \nabla \phi\rangle \right)e^{\frac{1}{a}\phi}d\upsilon=-\int_{M}|\nabla g|^{2}e^{\frac{1}{a}\phi}d\upsilon
\en
 which is a contradiction since
$g$ is the first eigenfunction of (\ref{1.5}) and cannot be a constant. Therefore, we have $\eta_{1}>\dfrac{\lambda_{1}}{n(a+1)}+b$. \hfill $\Box$ \vspace{0.2cm}
\\
\textbf{\emph{Proof of Theorem} \ref{Theorem 3}.} Let $w$ be an eigenfunction corresponding to the first
eigenvalue $q_{1}$ of the problem (\ref{1.10}), that is,
\be\no
\left\{ \begin{array}{ll}
L^2_{\phi}w=0 \ in \ M,   & \\
w=\frac{\partial^{2}u}{\partial\nu^{2}}-q_{1}\frac{\partial w}{\partial\nu}=0 \ on \ \partial M. & \\
\end{array} \right.
\en
Note that $w$ is not a constant since $w|_{\partial M}=0$. Let $\eta=\partial_{\nu}w|_{\partial M}$. Then $\eta\neq0$, otherwise, we would deduce from
\be\label{2.23}
w|_{\partial M}=\nabla w|_{\partial M}=\frac{\partial^{2}u}{\partial\nu^{2}}|_{\partial M}=0
\en
and (\ref{2.5}) that
\be\label{2.24}
L_{\phi}w|_{\partial M} = (\overline{L}_{\phi}w+(n-1)H_{\phi}w_{\nu}+\nabla^{2}w (\nu,\nu))|_{\partial M}=0.
\en
By divergence theorem, we have
\be\label{2.25}
0=\int_{\partial M}(L_{\phi}w\nabla L_{\phi}w\cdot\nu \ )e^{-\phi}dA=\int_{M}(|\nabla L_{\phi}w|^{2}+L_{\phi}wL^{2}_{\phi}w)e^{-\phi}d\upsilon,
\en
and so $L_{\phi}w=0$ on $M$. We also have
\be\label{2.26}
0=\int_{\partial M}(w\nabla w\cdot\nu) \ e^{-\phi}dA=\int_{M}(|\nabla w|^{2}+wL_{\phi}w)e^{-\phi}d\upsilon,
\en
that implies $w=0$. This is a contradiction.
From $w|\partial M=0$, one gets again from divergence theorem that
\be\label{2.27}
\int_{M}\langle\nabla w,\nabla L_{\phi}w\rangle e^{-\phi}d\upsilon=-\int_{M}wL^{2}_{\phi}we^{-\phi}d\upsilon=0,
\en
and therefore
\be\no\label{2.28}
\int_{\partial M}(L_{\phi}w\nabla w\cdot\nu) \ e^{-\phi}dA &=&\int_{M}(\langle\nabla L_{\phi}w,\nabla w \rangle+(L_{\phi}w)^{2})e^{-\phi}d\upsilon\\
&=&\int_{M}(L_{\phi}w)^{2}e^{-\phi}d\upsilon.
\en
Since that $\nabla^{2}w (\nu,\nu))|_{\partial M}=w_{\nu \nu}|_{\partial M}-(\nabla_{\nu}\nu)(w)|_{\partial M}=\dfrac{\partial^{2}w}{\partial\nu^{2}}|_{\partial M}-\langle\nabla w,\nabla_{\nu}\nu\rangle|_{\partial M}$,
by (\ref{2.5}) we get
\be\label{2.29}
L_{\phi}w|_{\partial M} = (n-1)H_{\phi}w_{\nu}+q_{1}w_{\nu},
\en
and together (\ref{2.28}) given us
\be\label{2.30}
q_{1}=\dfrac{\int_{M}(L_{\phi}w)^{2} \ d\mu-(n-1)\int_{\partial M}H_{\phi}\eta^{2} \ d\vartheta}{\int_{\partial M}\eta^{2} \ d\vartheta},
\en
which, combining with Reilly's formula and (\ref{2.12}), gives
\be\label{2.31}
q_{1}&\geq&\dfrac{1}{n(a+1)}\dfrac{\int_{M}(L_{\phi}w)^{2}d\mu}{\int_{\partial M}\eta^{2}d\vartheta}.
\en
On the other hand, we have from the variational characterization of $p_{1}$ (cf. (\ref{1.13}))
\be\label{2.32}
\dfrac{\int_{M}(L_{\phi} w)^{2} \ d\mu}{\int_{\partial M}\eta^{2} \ d\vartheta}\geq p_{1}.
\en
By (\ref{2.31}) and (\ref{2.32}), we get $q_{1}\geq \dfrac{p_{1}}{n(a+1)}$.
From [\cite{QX}, Theorem 1.3], we know
that equality holding if and only if $M$ is isometric to a ball in $\mathbb{R}^{n}$.
This proves the Theorem \ref{Theorem 3}.\hfill $\Box$ \vspace{0.2cm}
\\
\textbf{\emph{Proof of Theorem} \ref{Theorem 4}}. Let $f$ be an eigenfunction corresponding to the first eigen-
value $p_{1}$ of the problem (\ref{1.9}), that is
\be\label{2.33}
\left\{ \begin{array}{ll}
L^2_{\phi}f=0 \ in \ M,   & \\
f=L_{\phi}f-p_{1}\frac{\partial f}{\partial\nu}=0 \ on \ \partial M. & \\
\end{array} \right.
\en
Set $\eta=\frac{\partial f}{\partial \nu}|_{\partial M}$. Then
\be\label{2.34}
p_{1}=\frac{\int_{M} (L_{\phi}f)^{2}e^{-\phi}d\upsilon}{\int_{\partial M} \eta^{2}e^{-\phi}dA}.
\en
Substituting $f$ into Reilly$'$s formula and using (\ref{1.15}), we have
\be\no
\int_{M}((L_{\phi}f)^{2}-|\nabla^{2}f|^{2})e^{-\phi}d\upsilon\geq\int_{ M}\left(\dfrac{|\nabla\phi|^{2}}{na}-b\right)|\nabla f|^{2}+(n-1)c\int_{\partial M}\eta^{2}e^{-\phi}dA,
\en
By (\ref{2.12}) and (\ref{2.14}) we have
\be\no
\left(1-\dfrac{1}{n(a+1)}+\dfrac{b}{\lambda_{1}}\right)\int_{M}(L_{\phi}f)^{2}e^{-\phi}d\upsilon &\geq&(n-1)c\int_{\partial M}\eta^{2}e^{-\phi}dA,
\en
\be\label{2.35}
p_{1}\geq\dfrac{n(a+1)(n-1)c\lambda_{1}}{n(a+1)(\lambda_{1}+b)-\lambda_{1}}.
\en
If the equality sign holds in (\ref{2.35}), using the same arguments
as in the proof of \ref{Theorem 3} ([\cite{QX}, Theorem 1.3]), we conclude that  $M$
is isometric to a ball in $\mathbb{R}^{n}$ of radius $1/c$.\hfill $\Box$ \vspace{0.2cm}
\\
\textbf{\emph{Proof of Theorem} \ref{Theorem 5}}. From (\ref{1.20}), we only need to show that the first non-zero
eigenvalue $\xi_{1}$ of the problem \ref{1.11}) satisfies
\be\label{2.36}
\xi_{1}>\dfrac{cn(a+1)\lambda_{1}\mu_{1}}{n(a+1)(\mu_{1}+b)-\mu_{1}}.
\en
Let $f$ be an eigenfunction corresponding $\xi_{1}$:
\be\label{2.37}
\left\{ \begin{array}{ll}
L^2_{\phi}u=0 \ in \ M,   & \\
\frac{\partial f}{\partial \nu}=\frac{\partial (L_{\phi}f)}{\partial\nu}+\xi_{1} f=0 \ on \ \partial M. & \\
\end{array} \right.
\en
Let $z=f|_{\partial M}$; then $z\neq0$ and
\be\label{2.38}
\xi_{1}=\dfrac{\int_{M}(L_{\phi}f)^{2}d\mu}{\int_{\partial M}z^{2}d\vartheta}.
\en
Substituting $f$ into Reilly$'$s formula, we have
\be\no\label{2.39}
\int_{M}((L_{\phi}f)^{2}-|\nabla^{2}f|^{2})e^{-\phi}d\upsilon & = & \int_{ M}Ric_{\phi}(\nabla f,\nabla f)e^{-\phi}d\upsilon +\int_{\partial M}II(\overline{\nabla} z,\overline{\nabla} z)e^{-\phi}dA\\
\no& \geq & \int_{ M}\left(\dfrac{|\nabla \phi|^{2}|\nabla f|^{2}}{na}-b\right)|\nabla f|^{2}e^{-\phi}d\upsilon\\
&  &+c\int_{\partial M}|\overline{\nabla} z|^{2}e^{-\phi}dA.
\en
Since $\partial_{\nu}f|_{\partial M}=0$, we have
\be\label{2.40}
\int_{M}(L_{\phi}f)^{2}e^{-\phi}d\upsilon\geq\mu_{1}\int_{M}|\nabla f|^{2}e^{-\phi}d\upsilon.
\en
It follows from (\ref{2.37}) that $\int_{\partial M}z d\vartheta=0$. Indeed, by (\ref{2.37}) we have
\be
\no\xi_{1}\int_{\partial M}fe^{-\phi}dA & = & \int_{\partial M}\dfrac{\partial}{\partial \nu}(f-\overline{L}_{\phi} f)e^{-\phi}dA\\
\no& = &\int_{M}div((\nabla f-\nabla L_{\phi} f)e^{-\phi})d\upsilon\\
\no& = &\int_{M}(L_{\phi} f-L^{2}_{\phi} f)e^{-\phi}d\upsilon\\
\no& = &\int_{M}(L_{\phi} f)e^{-\phi}d\upsilon=0.
\en
So we have from the Poincar\'e inequality that
\be\label{2.41}
\int_{\partial M}|\overline{\nabla} z|^{2} \ e^{-\phi}dA\geq\lambda_{1}\int_{\partial M}z^{2} \ e^{-\phi}dA.
\en
Combining \ref{2.12}) and (\ref{2.38})-(\ref{2.41}), we get
\be\label{2.42}
\xi_{1}\geq\dfrac{cn(a+1)\lambda_{1}\mu_{1}}{n(a+1)(\mu_{1}+b)-\mu_{1}}.
\en
Let us show by contradiction that the equality in (\ref{2.42}) can’t occur. In fact, if
(\ref{2.42}) take equality sign, then we must have (\ref{2.11}) ocurring on $M$, that is,
\be\no
\nabla ^{2}f=\dfrac{\triangle f}{n}\langle,\rangle.
\en
Thus for a tangent vector field $X$ of $\partial M$, we have from by expression above and $\partial_{\nu}f|_{\partial M}=0$ that
\be\label{2.43}
0=\nabla ^{2}f(\nu,X)=X\nu f-(\nabla_{X}\nu)f=-\langle\nabla_{X}\nu,\overline{\nabla }z\rangle.
\en
In particular, we have
\be\no
II(\overline{\nabla }z,\overline{\nabla }z)=0.
\en
This is impossible since $II=cI$ and $z$ is not constant. This finishes the proof
Theorem \ref{Theorem 5}.\hfill $\Box$ \vspace{0.2cm}\\

The method used for the demonstration of the results is classic and has been widely used in articles in the bibliography. The B\'erard article \cite{BR} is a pioneering reference to the generalized Simons equation satisfied for the second fundamental form of an immersion in a Riemannian manifold. The Simons type inequalities used can be deduced from the B\'erard article.




\end{document}